\documentclass{ifacconf}
\usepackage{graphicx}      
\usepackage{natbib}        
\usepackage{amsmath}
\usepackage{amssymb}
\newtheorem{theorem}{Theorem}[section]
\newtheorem{lemma}[theorem]{Lemma}

\newtheorem{corollary}[theorem]{Corollary}

\newenvironment{example}[1][Example]{\begin{trivlist}
\item[\hskip \labelsep {\bfseries #1}]}{\end{trivlist}}
\newenvironment{remark}[1][Remark]{\begin{trivlist}
\item[\hskip \labelsep {\bfseries #1}]}{\end{trivlist}}

\begin{document}
\begin{frontmatter}

\title{An Adaptive Observer Design for Takagi-Sugeno type Nonlinear System\thanksref{footnoteinfo}} 

\thanks[footnoteinfo]{The work was supported by FP7 project Energy in Time (EiT) under the grant no. 608981}

\author[First,Second]{Krishnan Srinivasarengan} 
\author[First,Second]{Jos\'{e} Ragot} 
\author[First,Second]{Didier Maquin}
\author[First,Second]{Christophe Aubrun}
\address[First]{CRAN UMR 7039, CNRS, France}
\address[Second]{CRAN UMR 7039, Universit\'{e} de Lorraine, Vandoeuvre-l\`{e}s-Nancy, Cedex, France \\(e-mail: krishnan.srinivasarengan@univ-lorraine.fr, jose.ragot@univ-lorraine.fr, didier.maquin@univ-lorraine.fr, christophe.aubrun@univ-lorraine.fr)}

\begin{abstract}
Takagi-Sugeno (T-S) type of polytopic models have been used prominently in the literature to analyze nonlinear systems. With the sector nonlinearity approach, an exact representation of a nonlinear system within a sector could be obtained in a T-S form. Hence, a number of observer design strategies have been proposed for nonlinear systems using the T-S framework. In this work, a design strategy for adaptive observers is presented for a type of T-S systems with unknown parameters. The proposed approach improves upon the existing literature in two folds: reduce the computational burden and provide an algorithmic procedure that would seamlessly connect the state estimation and parameter estimation parts of the observer design. Lyapunov approach is used for the stability analysis and the design procedure. The results are illustrated on a simulation example.
\end{abstract}


\end{frontmatter}

\section{INTRODUCTION}
Nonlinear system analysis is an area of interest in the control community for a long time. One of the ways to do this is by representing a given nonlinear system in its equivalent form. Linear parameter varying systems are popular in this regard. Takagi-Sugeno (T-S) type of system model is quasi-LPV, in that the parameter is allowed to be one of the system variables, including, inputs, outputs, and states of the system. One of the approaches to obtain the T-S form from a given nonlinear system is the sector nonlinearity approach (\cite{ohtake2003fuzzy}). This approach results in a T-S form which is an exact representation of the original system within a sector. The sector restriction of the dynamics is a reasonable assumption as variables of a practical system are bounded.

Observer design for T-S models is an area that has attracted a lot of works in the last decade. A subclass of these observers considers unknown parameter estimation along with that of the state estimation (i.e., adaptive observers). Adaptive T-S observers could be classified into two types: those that consider an additive unknown parameter in the form of unknown inputs and those that consider multiplicative unknown parameters. In \cite{ichalal2009state}, an adaptive observer in the unknown input estimation form is proposed for a system with unmeasured premise variable. The common observer structures are either a Proportional Integral (PI) or Proportional Multiple Integral (PMI) (for the parameter estimation part), depending upon the order of derivatives of the parameter(s) that doesn't vanish (see, for instance, \cite{lendek2010adaptive}, \cite{ichalal2009simultaneous}). A typical application of this type of adaptive observers deals with actuator fault diagnosis (see \cite{marx2007design}).

For adaptive observers considering multiplicative unknown parameters, there have been different system model scenarios and observer structures. In \cite{lendek2010stability}, an adaptive observer is designed for the estimation of unmodeled dynamics in a T-S system. The observer structure proposed is inspired from that in \cite{cho_systematic_1997}. 
The present work uses the same starting point and improves upon the design procedure.

In \cite{delmotte2013fouling}, a T-S observer with parameter estimation was designed for a heat exchanger fouling detection problem. The fouling coefficients are modeled using an uncertainty factor similar to the unmodeled dynamics in \cite{lendek2010stability} and a polynomial fuzzy observer is designed. In \cite{bezzaoucha2013state}, a joint state and parameter estimation observer was proposed for T-S systems whose matrices depend on unknown parameters. In this strategy, the unknown parameters are rewritten using the sector nonlinearity approach through weighting functions and sector extremum values. To handle to situation arising out of the difference between the actual and estimated weighting functions (due to their dependence on unknown parameters/states), an $\mathcal{L}_2$ formulation is used. A first order structure is employed for the parameter estimation part of the observer. In \cite{srinivasarengan2016nonlinear}, the results were applied to a system model for a heat-exchanger-zones combination and the implementation challenges were discussed.

In this paper, an observer is proposed for systems of the model type that is a generalization of that considered in \cite{bezzaoucha2013state}. The observer structure is inspired by that proposed in \cite{cho_systematic_1997}. This observer structure is simpler from a computational perspective as compared to \cite{bezzaoucha2013state} and the procedure proposed would overcome the strong rank constraints on the transmission matrices as in \cite{lendek2010stability}. Lyapunov stability analysis is used to ensure the convergence of estimates under some persistency of excitation conditions.

The paper is organized as follows: in the following section, the system and observer model structures considered in this work and their motivations are described. The main results and the proof follow in Sec. \ref{sec_main_results}, which includes a discussion on the future extensions for the unknown premise variable case. The illustration of these results applied to a simulation example is given in Sec. \ref{sec_sim_example}. The Sec. \ref{sec_conclusions} concludes the paper by providing the future outlook.
\section{PRELIMINARIES} \label{sec_sys_models}
\subsection{Notations}
Takagi-Sugeno models are of the form,
\begin{align}
\dot{x}(t) &= \sum_{i=1}^r \mu_i(z(t))\lbrace A_i x(t) + B_i u(t)\rbrace \nonumber\\
y(t) &= Cx(t) \label{eq_ts_fund_model}
\end{align}
The paper uses the following notations and dimensions: 
\begin{align}
x \in \mathbb{R}^n, \ \ u \in \mathbb{R}^{n_u}, \ \ z \in \mathbb{R}^{n_p}, \ \ y \in \mathbb{R}^{n_y} \nonumber
\end{align}
The weighting functions $\mu_i(z)$ correspond to the weighted product of membership functions of each premise variable of the particular submodel $i$ (see for e.g., \cite{tanaka2004fuzzy} for more details). The number of submodels is $r = 2^{n_p}$. In this paper, the time factor $(t)$ is dropped in the expressions for simplicity in representation. The weighting functions satisfy the convex sum property, such that,
\begin{align}
\sum_i^r \mu_i(.) = 1 \qquad \text{and} \qquad 0\leq \mu_i(.) \leq 1, \  \forall i
\end{align}
\subsection{Motivation for system and observer model structures}
The T-S system model structure with unknown parameters $\theta_j$ considered in this work is of the form:
\begin{align}
\dot{x} &= \sum_{i=1}^{r} \mu_i(z) \lbrace(A_i+\sum_{j=1}^{n_\theta} \theta_j \bar{A}_{ij}) x + (B_i+\sum_{j=1}^{n_\theta} \theta_j \bar{B}_{ij}) u \nonumber  \\
& \qquad \qquad \qquad \qquad +(F_i + \sum_{j=1}^{n_\theta} \theta_j \bar{F}_{ij}) \rbrace \nonumber \\
y &= Cx \label{eq_sys_model}
\end{align}
While the unknown parameters are shown to affect all the matrices, the corresponding transmission matrices ($\bar{A}_{ij}, \bar{B}_{ij}, \bar{F}_{ij}$) could be zero, if a particular unknown parameter $\theta_j$ does not affect it. This model structure is a generalization of that in \cite{bezzaoucha2013state} and was motivated by a multitude of factors. In \cite{srinivasarengan2016nonlinear}, a practical scenario such problem structure could arise is illustrated.   In the case of the T-S models obtained through identification, fault diagnosis problem could be posed as a parameter estimation problem. Further, consider the nonlinear system of the form,
\begin{align}
\dot{x} &= Ax + \phi(x,u) + b f(x,u) \theta \nonumber \\
y &= Cx \label{eq_nonlin_sys_cho}
\end{align}
In \cite{cho_systematic_1997} and further in \cite{besancon_remarks_2000}, it has been shown that a state and parameter observer with asymptotically vanishing error is possible under certain conditions. Applying sector nonlinearity (SNL) on this system can transform it to the system model structure of the form \eqref{eq_sys_model}. The case when premise variables are not measured (i.e., $\hat{z}$ instead of $z$) could be handled through techniques discussed later. For this type of system, the objective is to design an observer with a Luenberger structure for the state estimation part and a time varying gain structure for the parameter estimation component. That is, 
\begin{align}
\dot{\hat{x}} &= \sum_{i=1}^{r} \mu_i(z) \lbrace(A_i+\sum_{j=1}^{n_\theta} \hat{\theta}_j \bar{A}_{ij}) \hat{x} + (B_i+\sum_{j=1}^{n_\theta} \hat{\theta}_j \bar{B}_{ij}) u \nonumber \\
& \qquad \qquad \qquad +(F_i + \sum_{j=1}^{n_\theta} \hat{\theta}_j \bar{F}_{ij}) + L_i(y-\hat{y})\rbrace \label{eq_obs_structure} \\
\dot{\hat{\theta}}_j &= f_j(\hat{\theta}_j, \hat{x}, y), \ \forall j =1,...,n_\theta \nonumber \\
\hat{y} &= C\hat{x} \nonumber
\end{align}
The choice of a model of this form is also an attempt to improve upon the overall computational burden on the joint state and parameter estimation for T-S systems. A time varying gain structure for the parameter estimation part allows for a more flexible dynamics, and at the same time avoids introducing complicated constraints as in \cite{bezzaoucha2013state}.

\section{RESULTS} \label{sec_main_results}
\subsection{Assumptions and Stability Analysis}
The observer design for the T-S system of the form \eqref{eq_sys_model} is given under the following assumptions:
\begin{enumerate}
	\item There exists a $\bar{\theta}$ such that $|\theta_j| < \bar{\theta}$, $\forall j$. This value is assumed to be known. However, the maximum value allowed by the design process could be determined as described later. \label{assumption_thetamax}
	\item All the submodels are sufficiently excited, illustrated by the variation in the weighting functions of each submodels, so that the system is under a persistence of excitation.
	\item $\dot{\theta} = 0$. The proof uses this condition to assume that the unknown parameters are constant. However, it is shown in the examples that this approach will work for slowly varying parameters as well.
\end{enumerate}

Computing the state estimation error between the system \eqref{eq_sys_model} with that of the Luenberger observer structure in \eqref{eq_obs_structure} as (with  $e_x \triangleq x-\hat{x}$ and $e_{\theta_j} \triangleq \theta_j -\hat{\theta}_j$),
\begin{align}
\dot{e}_x &= \sum_{i=1}^{r} \mu_i(z) \lbrace(A_i - L_iC)e_x + \sum_{j=1}^{n_\theta}(\theta_j \bar{A}_{ij} x - \hat{\theta}_j \bar{A}_{ij} \hat{x}) .\nonumber \\ 
& \qquad \qquad \qquad \qquad + (\bar{B}_{ij} u + \bar{F}_{ij})e_{\theta_j} \rbrace \nonumber
\end{align}
By adding and subtracting $ \sum_{j=1}^{n_\theta}\theta_j \bar{A}_{ij} \hat{x}$, the error dynamics becomes,
\begin{align}
\dot{e}_x &= \sum_{i=1}^{r} \mu_i(z) \lbrace (A_i-L_iC+\sum_{j=1}^{n_\theta} \theta_j \bar{A}_{ij}) e_x  \nonumber \\ & \qquad \qquad \qquad + (\sum_{j=1}^{n_\theta} (\bar{A}_{ij} \hat{x} +  \bar{B}_{ij} u + \bar{F}_{ij}) e_{\theta_j} \rbrace \nonumber
\end{align}
To analyze stability, consider the following Lyapunov function
\begin{align}
V &= e_x^T P e_x +  \sum_{j=1}^{n_\theta} e_{\theta_j}  \rho_j e_{\theta_j}
\end{align}
Its derivative is then given by
\begin{align}
\dot{V} &= \dot{e}_x^T P e_x + e_x^T P \dot{e}_x +  2 \sum_{j=1}^{n_\theta} \rho_j \dot{e}_{\theta_j} e_{\theta_j}\nonumber 
\end{align}
Considering,
{\small{
\begin{align*}
G_{ij} \triangleq P(A_i-L_i C+ \sum_{j=1}^{n_\theta} \theta_j \bar{A}_{ij}) + (A_i-L_i C + \sum_{j=1}^{n_\theta} \theta_j \bar{A}_{ij})^T P
\end{align*}
}}
and since $\dot{\theta} = 0$,  
\begin{align}
 2\sum_{j=1}^{n_\theta} \rho_j \dot{e}_{\theta_j} e_{\theta_j} = -2 \sum_{j=1}^{n_\theta} \rho_j \dot{\hat{\theta}}_j e_{\theta_j}
\end{align}
leads to,
\begin{align}
\dot{V} &= \sum_{i=1}^{r} \mu_i(z) \lbrace e_x^T G_{ij}  e_x \nonumber \\ 
& \ + 2 \sum_{j=1}^{n_\theta} e_{\theta_j} (\bar{A}_{ij} \hat{x} + \bar{B}_{ij} u + \bar{F}_{ij})^T P e_x\rbrace - 2  \sum_{j=1}^{n_\theta} \rho_j \dot{\hat{\theta}}_j e_{\theta_j} \label{eq_error_dyn}
\end{align}
To ensure $\dot{V} < 0$, each term on the right hand side (RHS) of \eqref{eq_error_dyn} is analyzed. The first contains a quadratic term with the unknown parameter in $G_{ij}$. Following the Assumption \ref{assumption_thetamax}, the result of robust stability analysis for error dynamics as in \cite{bergsten2001fuzzy} (and \cite{bergsten2000thau}) is applied. This translates to,
\begin{align}
P=P^T,\ P > 0, &\ Q=Q^T, \ Q > 0 \label{eq_conds_lmi1}\\
P(A_i - L_i C) +& (A_i - L_i C)^T P < -Q \label{eq_cond2_lmi}\\
n_\theta \bar{a} \bar{\theta} &\leq \frac{\lambda_{min}(Q)}{2 \lambda_{max}(P)} \label{eq_cond3_evals}
\end{align}
where, $\bar{a}$ is the maximum norm of all of $\bar{A}_{ij}, \ \forall i,j$.
The equivalent LMI condition for \eqref{eq_cond2_lmi} is given by:
\begin{align}
PA_i + A_i^T P- M_i C - C^T M_i^T < -Q \label{eq_conds_lmi2}
\end{align}
with the observer gain obtained as, $L_i = P^{-1} M_i$. The condition \eqref{eq_cond3_evals} will be satisfied if the following LMI condition is met, 
\begin{align}
\begin{bmatrix}
Q-\gamma I & P \\
P & I
\end{bmatrix} > 0 \label{eq_conds_lmi3}
\end{align} 
where, $\gamma = (n_\theta \bar{a} \bar{\theta})^2,\ \forall i,j$. The second and third terms in the RHS of \eqref{eq_error_dyn} relate to the coefficients of $e_\theta$ and one way to manage the Lyapunov function is to annihilate the coefficients of each error $e_{\theta_j}$,
\begin{align}
 \sum_{i=1}^{r}\mu_i(z) (\bar{A}_{ij} \hat{x} + \bar{B}_{ij} u +\bar{F}_{ij})^T P e_x - \rho_j \dot{\hat{\theta}}_j  = 0, \ \forall j \nonumber
\end{align}
This would lead to the condition,
\begin{align}
\dot{\hat{\theta}}_j &= \frac{1}{\rho_j} \sum_{i=1}^{r} \mu_i(z) \left( \bar{A}_{ij} \hat{x} + \bar{B}_{ij} u +\bar{F}_{ij}\right )^T P e_x, \ \ \forall\ j \label{eq_theta_estimate_e_x}
\end{align}
Since $e_x$ is not available, construction of $e_x$ from $e_y$ for the conditions specified is to be explored. This could be resolved in multiple ways, which are summarized as the main results. 
\subsection{Main Results}
The observer,
\begin{align}
\dot{\hat{x}} &= \sum_{i=1}^{r} \mu_i(z) \lbrace(A_i+\sum_{j=1}^{n_\theta} \hat{\theta}_j \bar{A}_{ij}) \hat{x} + (B_i+\sum_{j=1}^{n_\theta} \hat{\theta}_j \bar{B}_{ij}) u \nonumber \\
& \qquad \qquad \qquad +(F_i + \sum_{j=1}^{n_\theta} \hat{\theta}_j \bar{F}_{ij}) + L_i(y-\hat{y})\rbrace \nonumber \\
\dot{\hat{\theta}}_j &= \frac{1}{\rho_j} \sum_{i=1}^{r} \mu_i(z) \left( \bar{A}_{ij} \hat{x} + \bar{B}_{ij} u +\bar{F}_{ij}\right )^T P C^\dagger e_y, \ \ \forall\ j \nonumber \\
\hat{y} &= C\hat{x}  \label{eq_obs_model}
\end{align}
where $C^\dagger$ is the pseudo inverse of $C$, is an adaptive observer for the system \eqref{eq_sys_model} if one of the following theorems are satisfied, which resolve the problem in \eqref{eq_theta_estimate_e_x}.

\begin{theorem}(adapted from \cite{lendek2010stability})
The system \eqref{eq_obs_model} forms an observer for the system \eqref{eq_sys_model}, if
\begin{enumerate}
\item The conditions \eqref{eq_conds_lmi1}, \eqref{eq_conds_lmi2} and \eqref{eq_conds_lmi3} are satisfied
\item $N_i$ is of full column rank and 
\begin{align}
\text{rank}(C N_i) = \text{rank}(N_i), \ \forall i
\end{align}
where, $N_i$ is a pre-multiplying matrix that is common for all $\bar{A}_{ij}, \bar{B}_{ij}, \bar{F}_{ij} (\forall j)$ such that,
\begin{align}
\bar{A}_{ij} = N_i \tilde{A}_{ij}, \ \bar{B}_{ij} = N_i \tilde{B}_{ij}, \ \bar{F}_{ij} = N_i \tilde{F}_{ij}, \ \forall j
\end{align}
\end{enumerate}
Alternately this could be stated in the form of $\bar{A}_{ij}, \bar{B}_{ij}$, $\bar{F}_{ij} (\forall i,j)$ are of full column rank, and 
\begin{align}
\text{rank}(C\bar{A}_{ij}) &= \text{rank}(\bar{A}_{ij}) \nonumber \\
\text{rank}(C\bar{B}_{ij}) &= \text{rank}(\bar{B}_{ij}) \nonumber \\ 
\text{rank}(C\bar{F}_{ij}) &= \text{rank}(\bar{F}_{ij}) 
\end{align}
\end{theorem}
\paragraph*{Proof} Given the rank conditions, there exists, an $R_i$ such that, $R_i C = N_i^T P$ which could be used in \eqref{eq_theta_estimate_e_x} to give
\begin{align}
\dot{\hat{\theta}}_j &= \frac{1}{\rho_j} \sum_{i=1}^{r} \mu_i(z) \left( \tilde{A}_{ij} \hat{x} + \tilde{B}_{ij} u +\tilde{F}_{ij}\right )^T R_i  e_y, \ \ \forall\ j 
\end{align}
Considering $R_i = N_i^T P C^\dagger$, the parameter estimation part of the observer would become
\begin{align}
\dot{\hat{\theta}}_j &= \frac{1}{\rho_j} \sum_{i=1}^{r} \mu_i(z) \left( \bar{A}_{ij} \hat{x} + \bar{B}_{ij} u +\bar{F}_{ij}\right )^T P C^\dagger e_y, \ \ \forall\ j 
\end{align}
Hence the proof \qed
\begin{remark}[Remark 1]
This result is considerably restrictive due to the rank conditions on the system matrices of model chosen for this work. Apart from the structural constraints to be satisfied, there is no standard procedure that connects  choice of $R_i$ with that of the Lyapunov matrix $P$. To mitigate these problems, another approach is considered and given in the following theorem.
\end{remark}
\begin{theorem}
The system \eqref{eq_obs_model} forms an observer for the system \eqref{eq_sys_model}, if
\begin{enumerate}
\item The conditions \eqref{eq_conds_lmi1}, \eqref{eq_conds_lmi2} and \eqref{eq_conds_lmi3} are satisfied
\item For every $\theta_j$,
   \begin{align}
	\bar{A}_{ij}^T P H &= 0, \ \forall i \nonumber \\     
	\bar{B}_{ij}^T P H &= 0, \  \forall i \nonumber \\     
	\bar{F}_{ij}^T P H &= 0, \  \forall i \label{eq_lmEs}
	\end{align}
	where $H \triangleq I-C^\dagger C$ with $I$ being the identity matrix of appropriate dimensions.
\end{enumerate}
\end{theorem}
\paragraph*{Proof}
Consider the following lemma,
\begin{lemma}(\cite{rao1972generalized}) \label{lemma_generalized_solution_inverse}
Let $A \in \mathbb{R}^{m\times n}$ and $A^\dagger$ be any generalized inverse of $A$. Then a general solution of a consistent nonhomogeneous equation $A\mathbf{x}=\mathbf{y}$ is 
\begin{align}
\mathbf{x} = A^\dagger \mathbf{y} + H \omega
\end{align}
where $\omega$ is an arbitrary vector and $H=I-C^\dagger C$. The necessary and sufficient condition that $A\mathbf{x}=\mathbf{y}$ is consistent is,
\begin{align}
A A^\dagger \mathbf{y} = \mathbf{y}
\end{align}
\end{lemma}
The nonhomogeneous equation in the problem under focus is, $e_y = C e_x$, and based on this lemma,
\begin{align}
e_x = C^\dagger e_y + H\omega
\end{align}
for some arbitrary $\omega$. Applying this to \eqref{eq_theta_estimate_e_x} leads to,
\begin{align}
\hat{x}^T \bar{A}^T_{ij}  P e_x &=  \hat{x}^T \bar{A}^T_{ij} P (C^\dagger e_y + H\omega),  \ \forall i \nonumber \\ 
u^T \bar{B}^T_{ij} P e_x &= u^T \bar{B}^T_{ij} P (C^\dagger e_y + H\omega),  \ \forall i \nonumber \\
\bar{F}^T_{ij} P e_x &= \bar{F}^T_{ij} P (C^\dagger e_y + H\omega),  \ \forall i
\end{align}
The second term on the right hand side of the equation would lead to zero, if the conditions in \eqref{eq_lmEs} are satisfied. Hence, the proof. \qed
\begin{corollary}
If the value of $\bar{\theta}$ as in the Assumption \ref{assumption_thetamax} is not known, the maximum $\bar{\theta}$ allowed by a particular design could be obtained by rewriting the above results as an optimization problem considering $\gamma$ as an objective to be maximized. That is,
\begin{align}
\underset{P, Q, L_i}{\text{maximize}} \ \ \ \gamma \label{eq_opti_constraint}
\end{align}
$\forall i$, such that, the conditions \eqref{eq_conds_lmi1}, \eqref{eq_conds_lmi2} and \eqref{eq_conds_lmi3} are satisfied.
\end{corollary}
\begin{remark}[Remark 2] It is interesting to note the correlation between the structural constraints that arise out of this result with that of the the nonlinear adaptive observer form as proposed in \cite{besancon_remarks_2000}. In the reference, the state equations are split into those that are measured and unmeasured and the unknown parameter $\theta$ is allowed to appear only on the dynamics of the measured states. In \cite{cho_systematic_1997}, this constraint is given, for the system in \eqref{eq_nonlin_sys_cho}, as $b^TP$ should be in the space spanned by $C$. This is seen from the constraints in \eqref{eq_lmEs},
\begin{align}
\bar{A}_{ij}^T P  = (\bar{A}_{ij}^T PC^\dagger) C \Rightarrow \bar{A}_{ij}^T P \in \text{span}(C)
\end{align}
This approach also brings in a procedure where the design of $P$ is connected to the constraints on the system matrices that depend on $\theta$.
However, the equality constraint is restrictive and could be modified so as to minimize the sum of all the terms on the left hand side in the equality conditions, which leads to an optimization problem, 
\end{remark}
\begin{theorem} \label{thm_optimize}
The observer \eqref{eq_obs_model} for the system \eqref{eq_sys_model}, could be designed if the following optimization problem has a solution,
\begin{align}
\underset{P, M_i}{\text{min}} \ \sum_{j=1}^{n_\theta} \beta_j
\end{align}
under the constraints of \eqref{eq_conds_lmi1}, \eqref{eq_conds_lmi2}, \eqref{eq_conds_lmi3}. Here,
\begin{align}
\beta_j = \sum_{i=1}^r \Vert\bar{A}_{ij}^T P H\Vert+\Vert\bar{B}_{ij}^T P H\Vert+\Vert\bar{F}_{ij}^T P H\Vert
\end{align}
\end{theorem}
\begin{remark}[Remark 3]
The pseudo inverse of the output matrix, $C^\dagger$ could be computed through the Singular Value Decomposition (SVD). Given, $C = U\Sigma V^T$, then $C^\dagger = V \Sigma^\dagger U^T$.
\end{remark}
\begin{remark}[Remark 4]
It is to be noted that the parameter $\rho_j$ is not considered in the design process. It is manually tuned to improve the dynamics of the parameter estimation and the values would depend on the scale of the corresponding $\theta_j$.
\end{remark}
\begin{remark}[Remark 5]
If the T-S models are obtained through sector nonlinearity approach, there is a likelihood of the weighting functions would depend on estimated premise variables rather than measured. In this scenario, there are multiple ways to handle this. One popular approach is the use of Lipschitz condition (see for instance, \cite{lendek2010stability}), where the difference between the estimated and actual weighting functions are assumed to be bounded by a known value. In \cite{ichalal2016auxiliary}, a summary of different methods for observer design with unmeasured premise variables is provided. All these approaches are extendible to adaptive observers as well. Further, the authors propose a systematic approach using the immersion techniques to avoid unmeasureable premise variables in the resultant T-S model, if the nonlinearities have a polynomial structure. This approach shall be explored to extend the results in this paper.
\end{remark}
\vspace{-8pt}
\section{ILLUSTRATIVE EXAMPLE} \label{sec_sim_example}
The results of the proposed observer design is illustrated using the following example. 
\begin{example}[Example]
Consider a nonlinear system of the form
\begin{align}
\dot{x}_1 &= -0.7x_1^2 - x_2 + x_3 + (1-0.8x_1)\theta \nonumber\\
\dot{x}_2 &= -x_1x_3 - 2x_2 + (x_2+u) \theta \nonumber \\
\dot{x}_3 &= 0.5x_1 - 2x_3 + u \nonumber \\
y_1 &= x_1+x_2 \nonumber \\
y_2 &= x_2
\end{align}
With sector nonlinearity, this system could be transformed into a system of the form \eqref{eq_sys_model}. $z \triangleq x_1$ is the premise variable and is assumed to be in the sector of $(0,\ 2)$. The weighting functions are given by $\mu_1 = \frac{z_1}{2}$ and $\mu_2 = 1-\mu_1$. The following are the system matrices obtained after applying the sector nonlinearity approach.
\begin{align}
A_1 &= \begin{bmatrix}
-1.4 & -1 & 1 \\ 0 & -2 & -2 \\ 0.5 & 0 & -2
\end{bmatrix}, \ \ 
A_2 = \begin{bmatrix}
0 & -1 & 1 \\ 0 & -2 & 0 \\ 0.5 & 0 & -2
\end{bmatrix} \nonumber \\ 
\bar{A}_{11} &= \bar{A}_{21} = \begin{bmatrix}
-0.8 & 0 & 0 \\ 0 & 1 & 0 \\ 0 & 0 & 0
\end{bmatrix} \nonumber \\
B_1 &= B_2 = \begin{bmatrix}
0 \\ 0 \\ 1
\end{bmatrix},\bar{B}_{11} = \bar{B}_{12} = \begin{bmatrix}
0 \\ 1 \\ 0
\end{bmatrix}, \  
C = \begin{bmatrix}
1 & 1 & 0 \\ 0 & 1 & 0
\end{bmatrix} \nonumber
\end{align}
\end{example}

\subsection{Simulation Results}
The simulation of this example was carried out on Matlab with the \textit{Yalmip} modeling interface (\cite{lofberg2004yalmip}) and using the \textit{SeDuMi} solver (\cite{S98guide}). As could be seen, the conditions in \eqref{eq_lmEs} are satisfied for this example only if $P$ is constrained to be diagonal. By using Theorem \ref{thm_optimize}, the results show that this constraint is not required.

\begin{figure}[t]
\centering
\includegraphics[scale=0.5]{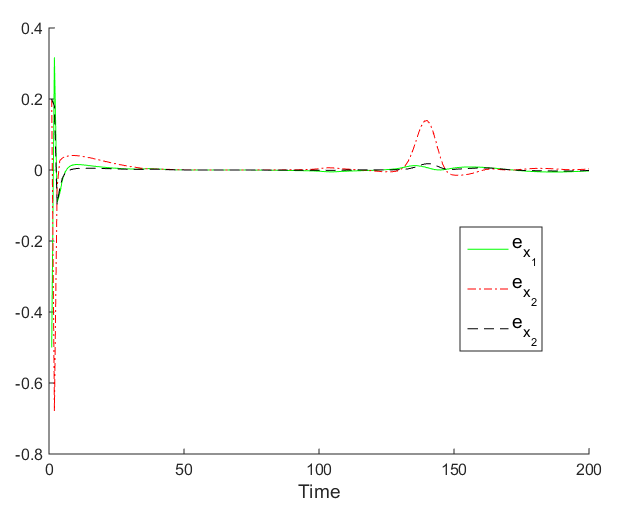}
\caption{State Errors' evolution over time}
\label{fig_state_est}
\end{figure}

\begin{figure}
\centering
\includegraphics[scale=0.55]{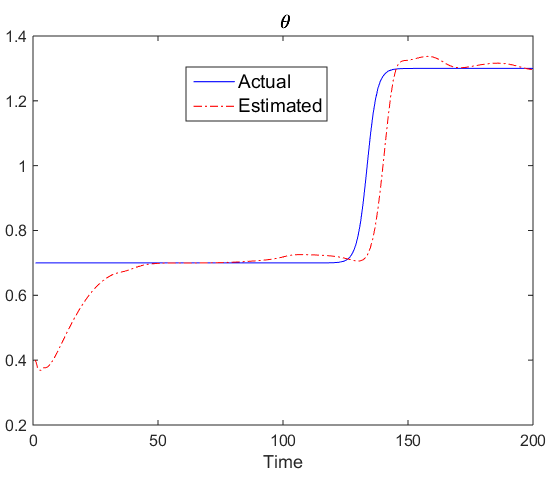}
\caption{Unknown parameter and its estimate}
\label{fig_para_est}
\end{figure}

\begin{figure}
\centering
\includegraphics[scale=0.58]{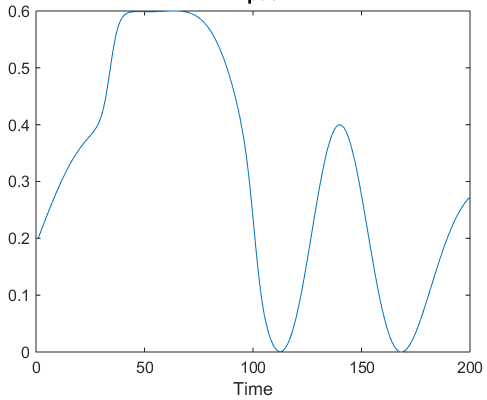}
\caption{Input used for the illustration}
\label{fig_input}
\end{figure}

\begin{figure}
\centering
\includegraphics[scale=0.55]{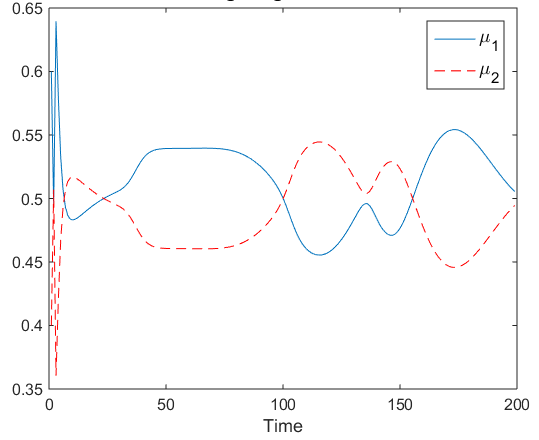}
\caption{Weighting function evolution for the simulation}
\label{fig_results_mu}
\end{figure}

The state estimation errors from this simulation are shown in the Fig. \ref{fig_state_est}. As could be seen, the estimation results are fairly accurate. The parameter estimation tracking is shown in the Fig. \ref{fig_para_est}, which shows a good tracking even when the unknown parameter changes ($\rho_1=1$ is used for this simulation). The input used for the simulation is shown in the Fig. \ref{fig_input} and the weighting functions in Fig. \ref{fig_results_mu} illustrate the sufficient excitation of both the submodels. The Lyapunov matrix obtained through the optimization problem was, 
\begin{align} \label{eq_P_example}
P = \begin{bmatrix}
1.169 & 0.657 & -4.7\times10^{-14} \\
0.657 & 1.153 & -3.3\times10^{-14} \\
 -4.7\times10^{-14} & -3.3\times10^{-14} & 1.365
\end{bmatrix}
\end{align}
and the value of the $\beta_1 = 1.31\times10^{-13}$. The structure obtained gives a key insight as discussed below.

\subsection{Discussion}
An observation about the connection between the structures of $C$ and $\bar{A}_{ij}, \bar{B}_{ij}, \bar{F}_{ij}$ with that of the Lyapunov matrix $P$ could be made. Consider the structure of $C$ of the form,
\begin{align}
C = \begin{bmatrix}
X_{n_y\times n_y} & & 0_{n_y\times(n_x-n_y)}
\end{bmatrix},
\end{align}
with $X \in \mathbb{R}^{n_y\times n_y}$ is a regular full rank matrix. This follows the following assumptions:
\begin{itemize}
\item The $C$ matrix is of full row rank. This is reasonable, for the redundant measurements, if exist, could be dropped.
\item There are some states that are not directly measured. As per the Remark 2, these are the states which do not have $\theta_j$ in their state equation.
\end{itemize}
Given $C$ is full column rank, $C^\dagger$ could be computed as 
\begin{align}
C^\dagger = C^T(CC^T)^{-1} = \begin{bmatrix} X^T \\ 0 \end{bmatrix} \begin{bmatrix} X X^T \end{bmatrix}^{-1}
\end{align}
This would lead to the matrix $H = I-C^\dagger C$,
\begin{align}
H &= I_{n_x} - \begin{bmatrix} X^T \\ 0 \end{bmatrix} \begin{bmatrix} X X^T \end{bmatrix}^{-1} \begin{bmatrix} X & 0 \end{bmatrix} = \begin{bmatrix} 0 & 0 \\ 0 & I_{n_x-n_y} \end{bmatrix} \nonumber
\end{align}
With this structure of $H$, some insights could be obtained for the structure of $P$. Consider, $P = \begin{bmatrix} P_1 & 0 \\ 0 & P_2 \end{bmatrix}$, 
with $P_1=P_1^T>0 \in \mathbb{R}^{n_y\times n_y}$ and $P_2>0$ is diagonal matrix of dimension $(n_x-n_y)$. If the transmission matrices ($\bar{A}_{ij}, \bar{B}_{ij}, \bar{F}_{ij}$) have a structure such that the unknown parameters affect only the measured states, then a Lyapunov matrix with the above structure would guarantee the equality conditions in \eqref{eq_lmEs}. However, there are no standard procedures to enforce such a structure on P. Fortunately, the optimization procedure in Theorem \ref{thm_optimize} facilitates this process without an explicit choice of the structure. This could be seen by the value of $P$ obtained as in \eqref{eq_P_example}, the procedure numerically tends to a structure of P with a full rank for the first $n_y = 2$ block and then towards a diagonal for the $n_x-n_y=1$ block. Hence it could be asserted that the algorithm facilitates an integration of connecting the stability requirements of the state estimation and the structural requirements of the parameter estimation.

\section{CONCLUDING REMARKS} \label{sec_conclusions}
In this paper, an adaptive observer for T-S models was proposed which has a time varying gain component for the parameter estimation part. The results on a simple example show that the performance is comparable to the existing literature results at a significantly less computational effort.

Moving forward, this work shall be extended to a case when the premise variables are not measured, but estimated using both known and upcoming approaches as discussed before. Further, the parameter estimation dynamics gain $\rho_j$ could be brought in as part of the design procedure to allow for a seamless design procedure. Taking cues from the proportional multiple integral (PMI) observers designed for unknown inputs, relaxing the condition of $\dot{\theta}=0$ to $\frac{d^n\theta}{dt^n} = 0$, for some $n>1$ could be explored.






%
%

\bibliography{./../../../../researchRef}
\end{document}